\documentclass[12pt]{article}

\usepackage{amsfonts}
\usepackage{amsmath}
\usepackage{amsthm}
\usepackage{graphicx}
\usepackage{ulem}

\newcommand{\RR}{\mathbb{R}}
\newcommand{\norm}[1]{\ensuremath{{\|#1\|}}}

\title{SDLS: a Matlab package for \\ solving conic
least-squares problems}
 
\author{Didier Henrion$^{1,2}$ \and J\'er\^ome Malick$^3$}

\begin{document}

\maketitle

\footnotetext[1]{LAAS-CNRS, University of Toulouse (France)}

\footnotetext[2]{Czech Technical University in Prague (Czech Republic)}

\footnotetext[3]{CNRS, Laboratoire Jean Kunztmann, University of Grenoble (France)}

\begin{abstract}
  This document is an introduction to the Matlab package SDLS (Semi-Definite Least-Squares) 
  for solving least-squares problems over convex symmetric cones.
  The package is shortly presented through the addressed problem,
  a sketch of the implemented algorithm, the syntax and calling sequences,
  a simple numerical example and some more advanced features. 
  The implemented method
  consists in solving the dual problem with a quasi-Newton algorithm.
  We note that SDLS is not the
  most competitive implementation of this algorithm:
  efficient, robust, commercial implementations are available
  (contact the authors). Our main goal with this Matlab SDLS package
  is to provide a simple, user-friendly software for solving and
  experimenting with semidefinite least-squares problems. Up to our knowledge,
  no such freeware exists at this date.
\end{abstract}

\section{General presentation}

SDLS is a Matlab freeware solving approximately
convex conic-constrained least-squares problems. Geometrically, such problems amount to
finding the best approximation of a point in the intersection of a convex symmetric cone 
with an affine subspace.
In mathematical terms, those problems can be cast as follows:
\begin{equation}\label{sdls}
\begin{array}{ll}
\min_x & \frac{1}{2}\|x-c\|^2 \\
\mathrm{s.t.} & Ax=b \\
& x \in K
\end{array}
\end{equation}
where $x \in {\mathbb R}^n$ has to be found,
$A \in {\mathbb R}^{m\times n}$, $b \in {\mathbb R}^m$, $c \in {\mathbb R}^n$
are given data and $K$ is a convex symmetric cone. 
The norm appearing in the objective function is  $\norm{x}=\sqrt{x^Tx}$,
the Euclidean norm associated with the standard inner product in $\RR^n$.
Note that there exists a unique solution to this optimization problem.

In practice, $K$ must be expressed as a direct product of linear,
quadratic and semidefinite cones.
It is always assumed that matrix $A$ has full row-rank, otherwise the problem
can be reformulated by eliminating redundant equality constraints.

\section{Algorithm}

SDLS is an implementation of the algorithm described in \cite{m04},
which consists in solving the dual to problem (\ref{sdls}). 
Specifically, the dual problem is (up to the sign and an additive constant)
\begin{equation}\label{dual}
\begin{array}{ll}
\min_y & f(y) = \frac{1}{2}\|p_K(c+A^Ty)\|^2-b^Ty \\
\end{array}
\end{equation}
where $p_K(z)$ denotes the orthogonal projection of vector $z$
onto the cone $K$. It is easy to prove that this dual problem 
is smooth and convex, so it can be solved by classical optimization 
algorithms. Among them, quasi-Newton algorithms are known to be efficient,
and SDLS uses a Matlab implementation
of the BFGS algorithm available in the Matlab freeware HANSO \cite{hanso}.
The other key numerical component of SDLS is eigenvalue decomposition for
symmetric matrices achieved by Matlab's built-in linear algebra function
$\tt eig$. We insist on the two following features:
\begin{itemize}
\item {\bf Simplicity:} the SDLS package consists of only 4 Matlab interpreted
 m-files calling one external package and a built-in linear algebra function.
\item {\bf Easy-to-use:} as explained in the next section, 
  the syntax follows standard patterns.
\end{itemize}

The implementation of the algorithm of \cite{m04} in SDLS 
is probably not the most efficient. It is not meant
to outperform neither the professional or commercial implementations of this method, 
nor other methods solving the same problem.
Our main goal is to provide a simple, user-friendly, free 
software for solving and experimenting with semidefinite least-squares problems.
Up to our knowledge, no such freeware exists at this date. 

\section{Syntax}\label{syntax}

SDLS has the same syntax as the widely used freeware SeDuMi
for solving linear problems over convex symmetric cones \cite{sedumi}.
Even though SeDuMi and SDLS share the same calling syntax, note that
SeDuMi is aimed at solving a different problem.
Namely, SeDuMi minimizes a linear objective function
$c^Tx$ while SDLS minimizes a quadratic objective function
as in~\eqref{sdls}.

The basic calling syntax is:
\begin{verbatim}
>> [x,y] = sdls(A,b,c,K)
\end{verbatim}

Input arguments {\tt A}, {\tt b}, {\tt c} are real-valued
matrices as in \eqref{sdls}. 
SDLS exploits data sparsity, so that {\tt A}, {\tt b} and {\tt c} can be
Matlab objects of class {\tt sparse}.
If the third input argument is empty
or not specified, then {\tt c} is set to zero.

Input argument {\tt K} is a structure describing the components of the cone:
\begin{itemize}
\item {\tt K.f} is the number of free, i.e. unrestricted components.
 E.g. if {\tt K.f=2} then {\tt x(1:2)} is unrestricted. These are
 always the first components in {\tt x}
\item {\tt K.l} is the dimension of the linear cone, i.e.
 the number of nonnegative components. E.g. if {\tt K.f=2} and {\tt K.l=8} then
 {\tt x(3:10)>=0}
\item {\tt K.q} lists the dimensions of the quadratic cones.
 E.g. if {\tt K.l=10} and {\tt K.q=[3 7]} then
 {\tt x(11)>=norm(x(12:13))} and {\tt x(14)>=norm(x(15:20))}.
 The corresponding components of {\tt x} immediately follow the {\tt K.l}
 nonnegative ones
\item {\tt K.s} lists the dimensions of the semidefinite cones.
 E.g. if {\tt K.l=20} and {\tt K.s=[4 3]}, then
 {\tt reshape(x(21:36),4,4)} and {\tt reshape(x(37:45),3,3)}
 are symmetric positive semidefinite matrices.
 These components are always the last entries in {\tt x}
\end{itemize}
Note that rotated quadratic cones ({\tt K.r})
and complex Hermitian components ({\tt K.xcomplex, K.scomplex, K.ycomplex})
are currently not supported. Please refer to SeDuMi's
documentation \cite{sedumi} for more details on the meaning
of the fields in structure {\tt K}.
If the fourth input argument is empty or not specified, then
{\tt K} is the linear cone, e.g. {\tt K.f=size(A,2)}.

Output argument {\tt x} (same dimension as {\tt c})
is the solution to SDLS problem (\ref{sdls}), whereas {\tt y} (same dimension
as {\tt b}) is the corresponding dual vector solving the smooth
problem (\ref{dual}).

\section{Installation}

SDLS 1.0 is developed for Matlab version 7.2 (release 2006).
It consists of a {\tt tar.gz} or {\tt zip} archive that can be
downloaded from
\begin{center}
{\tt www.laas.fr/$\sim$henrion/software/sdls}
\end{center}
Unpacking the archive 
creates a subdirectory {\tt sdls} that should be activated under
your Matlab environment, either with the inline command {\tt addpath} 
or with the {\tt Set Path} entry of the {\tt File} menu.

The package HANSO, available at
\begin{center}
{\tt www.cs.nyu.edu/overton/software/hanso}
\end{center}
must also be installed and activated under your Matlab environment.

If SDLS is correctly installed, it should generate the following
output:
\begin{verbatim}
>> x = sdls([1, 0],1)
x =
     1
     0
\end{verbatim}

\section{Example}

As a simple illustrative example of the use of SDLS
we consider the problem of calibrating correlation matrices.

Correlation matrices are symmetric positive semidefinite matrices
with ones along the diagonal.
They play an important role in probability and statistics
but also in combinatorial optimization.
Consider a symmetric matrix $C$ 
obtained from a correlation matrix after modifications of some entries,
for example due to measurement errors or post-processing.
Those modifications may alter the positive semidefiniteness of the matrix. 
We want to restore basic properties (positive semidefiniteness and ones along the
diagonal) by computing the correlation matrix nearest to $C$ in the
Euclidean sense. This problem can written in the form~\eqref{sdls}.

For illustration, consider the problem of computing a nearest
correlation matrix of size~3. Define
\begin{verbatim}
>> A=sparse(3,3); A(1,1)=1; A(2,5)=1; A(3,9)=1;
>> b=ones(3,1); % ones along diagonal
>> C=[1 1/2 1;1/2 1 1/4;1 1/4 1]
C =
    1.0000    0.5000    1.0000
    0.5000    1.0000    0.2500
    1.0000    0.2500    1.0000
>> c=C(:);
>> K=[]; K.s=3
\end{verbatim}
Matrix {\tt C} was obtained by adding {\tt 1/2} to entries {\tt (1,3)}
and {\tt (3,1)} of a correlation matrix. Modifying these 2 entries
destroyed the positive semidefiniteness:
\begin{verbatim}
>> min(eig(C))
ans =
   -0.0349
\end{verbatim}
Then we run the calibration process using SDLS:
\begin{verbatim}
>> [x,y] = sdls(A,b,c,K);
>> X = reshape(x,K.s,K.s)
X =
    1.0000    0.4910    0.9684
    0.4910    1.0000    0.2582
    0.9684    0.2582    1.0000
\end{verbatim}
This way, we do not recover the original correlation matrix, but we get
the correlation matrix nearest to {\tt C} in the Euclidean sense.
Indeed we observe that the diagonal ones remain whereas positive
semidefiniteness is restored:
\begin{verbatim}
>> min(eig(X))
ans =
   3.0097e-16
\end{verbatim}

\section{Advanced use}

\subsection{Input arguments}

Optionally, tuning parameters can be specified as a fifth
input argument:
\begin{verbatim}
>> [x,y] = sdls(A,b,c,K,pars)
\end{verbatim}
Input argument {\tt pars} is a structure containing the following fields:
\begin{itemize}
\item {\tt pars.eps}: relative accuracy (default {\tt 1e-6})
\item {\tt pars.fid}: screen output (0 for no output, default 1)
\item {\tt pars.maxit}: maximum number of iterations for the
BFGS solver (default {\tt 200})
\item {\tt pars.reg}: regularization parameter (default {\tt 0})
\item {\tt pars.scaling}: scaling of the linear equation
  ({\tt false} for no scaling, default {\tt true})

\end{itemize}

The stopping rule of the algorithm is {\tt
norm(A*x-b) <= max(1,norm(b))*pars.eps}
where {\tt pars.eps} is the expected relative accuracy.
Since BFGS is essentially a first-order algorithm, {\tt pars.eps}
should not be set too small. Typical values are between {\tt 1e-4}
and {\tt 1e-8}, depending on the problem dimensions and scaling.

The regularization parameter {\tt pars.reg} should be used if one suspects
that the qualification constraint does not hold. That is, the affine subspace
may be tangent to the cone, or equivalently there may be no point satisfying the
affine constraint in the interior of the cone.
A typical value is {\tt pars.reg=1e-6}.

When {\tt pars.scaling=true}, the scaling consists simply in dividing
input data {\tt A} and {\tt b} by {\tt max(1,norm(b))}.
Otherwise, the input data are left unchanged.

\subsection{Output arguments}

An additional third output argument:
\begin{verbatim}
>> [x,y,info] = sdls(A,b,c,K)
\end{verbatim}
can be retrieved. Argument {\tt info} contains information on
the objective function of dual problem
(\ref{dual}), as provided by the BFGS solver:
\begin{itemize}
\item {\tt info.f}: function value $f(y)$ at the optimum
\item {\tt info.g}: gradient of $f(y)$ at the optimum 
\item {\tt info.H}: approximation of inverse Hessian of $f(y)$ at the optimum
\item {\tt info.time}: CPU time
\end{itemize}

The norm of the residual {\tt norm(A*x-b)} is equal to 
the norm of the gradient {\tt norm(info.g)}. It should be less than
the absolute accuracy {\tt max(1,norm(b))*pars.eps} (see above)
and it can be used as an estimate
of the quality of the solution produced by SDLS.


\end{document}